\documentclass[12pt,a4paper]{article}

\newtheorem{theorem}{Theorem}[section]
\newtheorem{proposition}[theorem]{Proposition}
\newtheorem{lemma}[theorem]{Lemma}
\newtheorem{corollary}[theorem]{Corollary}

\usepackage{hyperref}
\usepackage{amsmath}
\usepackage{amssymb}
\usepackage{array}
\usepackage[english]{babel}
\usepackage{calc}
\usepackage{enumerate}
\usepackage[latin1]{inputenc}
\usepackage{latexsym}

\begin{document}
\title{Properties of derivations on some convolution algebras}
\author{Thomas Vils Pedersen}

%\date
\maketitle


\begin{abstract}
\noindent
For all the convolution algebras $L^1[0,1),\ L^1_{\text{loc}}$ and $A(\omega)=\bigcap_n L^1(\omega_n)$, the derivations are of the form $D_{\mu} f=Xf*\mu$ for suitable measures $\mu$, where $(Xf)(t)=tf(t)$. We describe the (weakly) compact as well as the (weakly) Montel derivations on these algebras in terms of properties of the measure $\mu$. Moreover, for all these algebras we show that the extension of $D_{\mu}$ to a natural dual space is weak-star continuous.
\end{abstract}

\section{Introduction}
\label{sec:intro}

The aim of this paper is to study various properties of derivations on some convolution Banach and Fr\'echet algebras. A starting point for the paper is the characterisation in \cite[Theorem~4.1]{De-Gh-Gr} of (weak) compactness of derivations on the weighted convolution Banach algebras $L^1(\omega)$. Other inspirations include the recent papers \cite{Ch-He:Trans}, \cite{Ch-He} and  \cite{Pe:Der-prop} on (weak) compactness and weak-star continuity of derivations from some Banach algebras to their dual spaces. 

The algebras that we will consider are $L^1[0,1),\ L^1_{\text{loc}}$ and $A(\omega)=\bigcap_n L^1(\omega_n)$ (see the relevant sections for the definitions). These are all convolution algebras on $[0,1)$ or $\mathbb R^+=[0,\infty)$ with the usual convolution product
$$(f\ast g)(t)=\int_0^t f(s)g(t-s)\,ds$$
for $f$ and $g$ in the algebra in question and $t\in[0,1)$ or $t\in\mathbb R^+$. On each of these algebras $\mathcal B$ all derivations are continuous and are of the form $D_{\mu} f=Xf*\mu\ (f\in\mathcal B)$ for $\mu$ in a suitable class of measures, where $X$ is the operator defined by $(Xf)(t)=tf(t)$ for $t\in[0,1)$ or $t\in\mathbb R^+$ and $f\in\mathcal B$. For derivations on $L^1[0,1)$ we obtain a characterisation of (weak) compactness in terms of the measure $\mu$ similar to the one for $L^1(\omega)$ in \cite{De-Gh-Gr}. For the Fr\'echet algebras $L^1_{\text{loc}}$ and $A(\omega)$ we will see that there are no non-zero weakly compact derivations, and in these cases the following seems to be a more useful notion: A linear operator between two Fr\'echet spaces is called {\it (weak) Montel} (see \cite{Bo-Li}) if it maps bounded sets to (weakly) relatively compact sets. On Banach spaces this notion agrees with the one of (weakly) compact operators, and generally a (weakly) compact operator is (weakly) Montel. For the Fr\'echet algebras $L^1_{\text{loc}}$ and $A(\omega)$ we give characterisations of (weak) Montel derivations similar to the characterisations of (weak) compactness for $L^1(\omega)$ and $L^1[0,1)$.

We also study weak-star continuity of the extension $\overline{D}_{\mu}$ of $D_{\mu}$ to a natural dual space containing the algebra in question. In all cases we prove weak-star continuity of $\overline{D}_{\mu}$ by showing that $\overline{D}_{\mu}$ is the adjoint of the continuous linear operator $T_{\mu}$ defined by
$$(T_{\mu} h)(t)=t\int_{\cdots}h(t+s)\,d\mu(s),$$
for $h$ belonging to the predual and $t\in[0,1)$ or $t\in\mathbb R^+$, and where the integrals are over $[0,1-t)$ or $\mathbb R^+$.

\section{Derivations on $L^1[0,1)$}
\label{sec:vol}

Let $L^1[0,1)$ be the Volterra algebra of (equivalence classes of) integrable functions $f$ on $[0,1)$ with convolution product and the norm 
$\|f\|=\int_0^1|f(t)|\,dt$. Similarly, $M[0,1)$ denotes the Banach algebra of finite, complex Borel measures on $[0,1)$. Also, we let $C_0[0,1)$ be the space of continuous functions $h$ on $[0,1]$ with $h(1)=0$. It is well known that
$\langle h,\mu\rangle=\int_{[0,1)}h(t)\,d\mu(t)$ for $h\in C_0[0,1)$ and $\mu\in M[0,1)$ identifies $M[0,1)$ isometrically isomorphically with the dual space of $C_0[0,1)$. 

The continuous derivations on $L^1[0,1)$ were described as follows by Kamowitz and Scheiberg (\cite[Theorem~2]{Ka-Sc}): Let $\mu$ be a measure on $[0,1)$ for which $t|\mu|([0,1-t))$ is bounded as $t\to0_+$. Then
$$D_{\mu} f=Xf*\mu\qquad(f\in L^1[0,1))$$
defines a continuous derivation on $L^1[0,1)$, and conversely every continuous derivation on $L^1[0,1)$ is of this form. Subsequently, Jewell and Sinclair (\cite{Je-Si}) proved that derivations on $L^1[0,1)$ are automatically continuous.

We first show that the derivations $D_{\mu}$ extend to weak-star continuous derivations on the measure algebra $M[0,1)$.

\begin{proposition}
\label{pr:mol}
Let $\mu$ be a measure on $[0,1)$ for which $t|\mu|([0,1-t))$ is bounded as $t\to0_+$. Then
$$\overline{D}_{\mu}\nu=X\nu*\mu\qquad(\nu\in M[0,1))$$
extends $D_{\mu}$ to a weak-star continuous, bounded derivation $\overline{D}_{\mu}$ on $M[0,1)$.
\end{proposition}

\noindent{\bf Proof}\quad
Let $\nu\in M[0,1)$ and $h\in C_0[0,1)$. Then
\begin{align*}
\left|\int_{[0,1)}h(t)\,d(X\nu*\mu)(t)\right|  
&= \left|\int_{[0,1)}\int_{[0,1-t)}h(t+s)\,d\mu(s)\,t\,d\nu(t)\right|\\  
&\le \int_{[0,1)}\int_{[0,1-t)}|h(t+s)|\,d|\mu|(s)\,t\,d|\nu|(s)\\
&\le \|h\|\int_{[0,1)}t|\mu|([0,1-t))\,d|\nu|(t)<\infty.
\end{align*}
This shows that $X\nu*\mu\in M[0,1)$ and that $\overline{D}_{\mu}$ is a bounded, linear operator on $M[0,1)$. A direct calculation shows that $X$ is a derivation on $M[0,1)$ and the same thus holds for $\overline{D}_{\mu}$.
\medskip

Let $h\in C_0[0,1)$ and let
$$(T_{\mu} h)(t)=t\int_{[0,1-t)}h(t+s)\,d\mu(s)\qquad(t\in[0,1]).$$
The continuity of $T_{\mu} h$ is relatively standard (see, for instance, \cite[Theorem~3.3.15]{Da:Book}), but there is a slight complication because we are integrating over $[0,1-t)$. For $0\le t,t_0\le 1$ we have
\begin{gather*}
|(T_{\mu} h)(t)-(T_{\mu} h)(t_0)
= \left|t\int_{[0,1-t)}h(t+s)\,d\mu(s)-t_0\int_{[0,1-t_0)}h(t_0+s)\,d\mu(s)\right|\\
\le |t-t_0|\int_{[0,1-t)}|h(t+s)|\,d|\mu|(s)
+t_0\int_{[0,1-t)}|h(t+s)-h(t_0+s)|\,d|\mu|(s)\\
+t_0\left|\int_{[0,1-t)}h(t_0+s)\,d\mu(s)-\int_{[0,1-t_0)}h(t_0+s)\,d\mu(s)\right|.
\end{gather*}
The first two terms tend to zero as $t\to t_0$. For $t\to t_0$ with $t\ge t_0$ the third term 
$t_0|\int_{[1-t,1-t_0)}h(t_0+s)\,d\mu(s)|\to0$, since $\bigcap_{t>t_0}[1-t,1-t_0)=\emptyset$. Similarly, for $t\to t_0$ with $t\le t_0$ the third term 
$t_0|\int_{[1-t_0,1-t)}h(t_0+s)\,d\mu(s)|\to t_0|h(1)\mu(\{1-t_0\})|=0$, since $h(1)=0$.
Hence $T_{\mu} h\in C_0[0,1)$ and it follows that $T_{\mu}$ is a continuous linear operator on $C_0[0,1)$. Moreover, for $\nu\in M[0,1)$ and $h\in C_0[0,1)$ we have
\begin{align*}
\langle h,\overline{D}_{\mu}\nu\rangle
&= \int_{[0,1)}h(t)\,d(X\nu*\mu)(t)\\  
&= \int_{[0,1)}\int_{[0,1-t)}h(t+s)\,d\mu(s)\,t\,d\nu(t)\\  
&= \langle T_{\mu} h,\nu\rangle.
\end{align*}
Hence $\overline{D}_{\mu}=T_{\mu}^*$ and in particular $\overline{D}_{\mu}$ is weak-star continuous.
{\nopagebreak\hfill\raggedleft$\Box$\bigskip}

The following characterisation of (weakly) compact derivations on $L^1[0,1)$ and its proof are strongly inspired by \cite[Theorem~4.1]{De-Gh-Gr}.

\begin{theorem}
\label{th:cptvol}
Let $\mu$ be a measure on $[0,1)$ for which $t|\mu|([0,1-t))$ is bounded as $t\to0_+$. Then the following conditions are equivalent:
\begin{enumerate}[(a)]
\item 
$D_{\mu}$ is a compact derivation on $L^1[0,1)$.
\item 
$D_{\mu}$ is a weakly compact derivation on $L^1[0,1)$.
\item
$\mu$ is absolutely continuous and $t|\mu|([0,1-t))\to0$ as $t\to0_+$.
\item 
$\overline{D}_{\mu}$ is a compact derivation on $M[0,1)$.
\item 
$\overline{D}_{\mu}$ is a weakly compact derivation on $M[0,1)$.
\end{enumerate}
\end{theorem}

\noindent{\bf Proof}\quad
The implications (d)$\Rightarrow$(a)$\Rightarrow$(b) and (d)$\Rightarrow$(e)$\Rightarrow$(b) are obvious.

(b)$\Rightarrow$(c):\quad
Let $t>0$, let $\delta_t$ be the Dirac point measure at $t$ and let $(e_k)$ be a bounded approximate identity for $L^1[0,1)$. Since $D_{\mu}$ is weakly compact there exist a subsequence $(e_{k_j})$ and $f\in L^1[0,1)$ such that
$$D_{\mu}(\delta_t*e_{k_j})\to f\qquad\text{weakly in $L^1[0,1)$ as $j\to\infty$}.$$
Let $g\in L^1[0,1)$. Since $L^1[0,1)^*$ is a $L^1[0,1)$ module we have $D_{\mu}(\delta_t*e_{k_j})*g\to f*g$ weakly in $L^1[0,1)$ as $j\to\infty$. Also,
$$D_{\mu}(\delta_t*e_{k_j})*g=\overline{D}_{\mu}(\delta_t)*e_{k_j}*g+\delta_t*D_{\mu}(e_{k_j})*g
\to\overline{D}_{\mu}(\delta_t)*g$$
in $L^1[0,1)$ as $j\to\infty$, since 
$D_{\mu}(e_k)*g=D_{\mu}(e_k*g)-e_k*D_{\mu}(g)\to D_{\mu}(g)-D_{\mu}(g)=0$ in $L^1[0,1)$ as $k\to\infty$. Hence $\overline{D}_{\mu}(\delta_t)*g=f*g$ for all $g\in L^1[0,1)$, so we deduce that $t\delta_t*\mu=\overline{D}_{\mu}(\delta_t)=f\in L^1[0,1)$. Since this holds for all $t>0$, it follows that $\mu$ is absolutely continuous.

Since the set $\{\delta_t*e_k:t>0\text{ and }k\in\mathbb N\}$ is bounded, it follows that $\{D_{\mu}(\delta_t*e_k):t>0\text{ and }k\in\mathbb N\}$ is weakly relatively compact in $L^1[0,1)$. Also, for $t>0$ we saw above that $t\delta_t*\mu=\overline{D}_{\mu}(\delta_t)$ is a weak cluster point of the sequence $(D_{\mu}(\delta_t*e_k))$, so we deduce that $\{t\delta_t*\mu:t>0\}$ is weakly relatively compact in $L^1[0,1)$. From the Dunford-Pettis characterisation 
(\cite[Theorem~4.21.2]{Ed:FA}) of weakly relatively compact subsets of $L^1[0,1)$
(or, as in \cite{De-Gh-Gr}, the Dieudonn\'e-Grothendieck characterisation 
(\cite[Theorem~4.22.1(4)]{Ed:FA}) of weakly relatively compact subsets of $M[0,1)$) it then follows that $\{t\delta_t*|\mu|:t>0\}$ is weakly relatively compact in $L^1[0,1)$. Let $(t_i)$ be any net in $(0,1)$ with $t_i\to0$. Then there exist a subnet $(t_{i_j})$ and $f\in L^1[0,1)$ such that
$t_{i_j}\delta_{t_{i_j}}*|\mu|\to f$ weakly in $L^1[0,1)$. Let $a>0$. Clearly $t_{i_j}\delta_{t_{i_j}}*\delta_a*|\mu|\to\delta_a*f$ weakly in $L^1[0,1)$, but $t_{i_j}\delta_{t_{i_j}}\to0$ in $M[0,1)$ and $\delta_a*|\mu|\in L^1[0,1)$, so $t_{i_j}\delta_{t_{i_j}}*\delta_a*|\mu|\to0$ in $L^1[0,1)$, 
% Moreover, $t_{i_j}\delta_{t_{i_j}}\to0$ weak-star in $M[0,1)$. We have 
% $\delta_a*|\mu|\in L^1[0,1)$ 
% and for $h\in C_0[0,1)$ an easy calculation shows that
% $$\langle h,t_{i_j}\delta_{t_{i_j}}*\delta_a*|\mu|\rangle
% =\langle h*(\delta_a*|\mu|),t_{i_j}\delta_{t_{i_j}}\rangle\to0$$
% since $h*(\delta_a*|\mu|)\in C_0[0,1)$. Hence $t_{i_j}\delta_{t_{i_j}}*\delta_a*|\mu|\to0$ 
% weak-star in $M[0,1)$,
and we deduce that $\delta_a*f=0$. Since this holds for all $a>0$ we have $f=0$, so we conclude that $t\delta_t*|\mu|\to 0$ weakly in $L^1[0,1)$ as $t\to0_+$. The constant function with value $1$ belongs to $L^{\infty}[0,1)=L^1[0,1)^*$, so we have 
$t|\mu|([0,1-t))=\langle t\delta_t*|\mu|,1\rangle\to 0$ as $t\to0_+$.

(c)$\Rightarrow$(d):\quad
We first prove that
$$E=\{t\delta_t*\mu:t\in[0,1)\}=\{\overline{D}_{\mu}(\delta_t):t\in[0,1)\}$$
is compact in $M[0,1)$. Let $(t_n\delta_{t_n}*\mu)$ be a sequence i $E$. We may assume that there exists $t_0\in[0,1]$ such that $t_n\to t_0$ as $n\to\infty$. Assume that $t_0=0$. Since $\|\overline{D}_{\mu}(\delta_t)\|=t|\mu|([0,1-t))$ we have 
$$t_n\delta_{t_n}*\mu=\overline{D}_{\mu}(\delta_{t_n})\to0\in E$$
as $n\to\infty$. Now assume that $0<t_0\le1$. Choose $t$ with $0<t<t_0$. Since $t\delta_t*\mu=\overline{D}_{\mu}(\delta_t)\in M[0,1)$ and since $\mu$ is absolutely continuous, we have $\delta_t*\mu\in L^1[0,1)$. Since $(\delta_s)$ is strongly continuous on $L^1[0,1)$, it follows that
$$t_n\delta_{t_n}*\mu=t_n\delta_{t_n-t}\delta_t*\mu
\to t_0\delta_{t_0-t}\delta_t*\mu=t_0\delta_{t_0}*\mu\in E$$
as $n\to\infty$ (with $\delta_1*\mu=0\in E$ in case $t_0=1$). Consequently $E$ is compact, so by Mazur's theorem (\cite[Theorem~VI.4.8]{Co}) the closed convex hull $\overline{\text{co}}(E)$ is compact. Let $F\subseteq M[0,1)$ consist of those finite point measures
$$\nu=\sum_{k=1}^K\alpha_k\delta_{t_k}\qquad
(K\in\mathbb N,\ 0\le t_1<t_2<\ldots<t_K<1)$$
for which $\|\nu\|=\sum_{k=1}^K|\alpha_k|\le1$. 
Then $\overline{D}_{\mu}(F)\subseteq\text{co}(E)$. Let $\nu\in M[0,1)$ with $\|\nu\|\le1$. As in the proof of \cite[Theorem~4.1]{De-Gh-Gr} there exists a net $(\nu_i)$ in $F$ with $\nu_i\to\nu$ strongly in $M[0,1)$. The net $(\overline{D}_{\mu}(\nu_i))$ belongs to the compact set $\overline{\text{co}}(E)$ and thus have a convergent subnet $(\overline{D}_{\mu}(\nu_{i_j}))$ with limit $\rho\in\overline{\text{co}}(E)$. Also, for $f\in L^1[0,1)$ we have
$$\overline{D}_{\mu}(\nu_{i_j})*f=D_{\mu}(\nu_{i_j}*f)-\nu_{i_j}*D_{\mu}(f)
\to D_{\mu}(\nu*f)-\nu*D_{\mu}(f)=\overline{D}_{\mu}(\nu)*f,$$
so we deduce that $\overline{D}_{\mu}(\nu)=\rho\in\overline{\text{co}}(E)$. Hence $\overline{D}_{\mu}$ is compact.
{\nopagebreak\hfill\raggedleft$\Box$\bigskip}

\section{Derivations on $L^1_{\text{loc}}$}
\label{sec:loc}

We denote by $L^1_{\text{loc}}$ the space of locally integrable functions on $\mathbb R^+$ and by $M_{\text{loc}}$ the space of Radon measures on $\mathbb R^+$, that is, locally finite, complex Borel measures on $\mathbb R^+$. For $n\in\mathbb N$ we define the restriction map $R_n:M_{\text{loc}}\to M[0,n)$ and the inclusion map $S_n:M[0,n)\to M_{\text{loc}}$ in the obvious ways. Equipped with the seminorms $\mu\mapsto\|R_n\mu\|\ (\mu\in M_{\text{loc}})$ for $n\in\mathbb N$ it is well known that $L^1_{\text{loc}}$ and $M_{\text{loc}}$ become Fr\'echet convolution algebras on $\mathbb R^+$. These algebras can also be regarded as the projective limits of the spaces $L^1[0,n)$ respectively $M[0,n)$.

The multipliers and derivations on $L^1_{\text{loc}}$ were described in 
\cite[Theorem~2.14 and Theorem~3.1]{Gh-Mc}: For $\mu\in M_{\text{loc}}$ the linear map $M_{\mu}f=\mu*f\ (f\in L^1_{\text{loc}})$ defines a continuous multiplier on $L^1_{\text{loc}}$ and conversely every multiplier on $L^1_{\text{loc}}$ is of this form. Similarly, for $\mu\in M_{\text{loc}}$ the linear map $D_{\mu}f=Xf*\mu\ (f\in L^1_{\text{loc}})$ defines a continuous derivation on $L^1_{\text{loc}}$ and conversely every derivation on $L^1_{\text{loc}}$ is of this form. (In particular, multipliers and derivations on $L^1_{\text{loc}}$ are automatically continuous.) Moreover, $\overline{D}_{\mu}\nu=X\nu*\mu\ (\nu\in M_{\text{loc}})$ extends $D_{\mu}$ to a continuous derivation on $M_{\text{loc}}$.

Let $C_c$ be the space of compactly supported, continuous functions on $\mathbb R^+$. We regard $C_c$ as the inductive limit of the spaces $C_0[0,n)$ and equip it with the corresponding inductive limit topology. It follows as in the proof of \cite[Proposition~3.3]{Pe:Fr} (see also \cite{Gr-hom}) that 
$$\langle h,\mu\rangle=\int_{\mathbb R^+}h(t)\,d\mu(t)\qquad 
(h\in C_c,\mu\in M_{\text{loc}})$$
identifies $M_{\text{loc}}$ with the dual space of $C_c$.

\begin{proposition}
\label{pr:locwks}
Let $\mu\in M_{\text{loc}}$. Then the derivation $\overline{D}_{\mu}$ is  weak-star continuous on $M_{\text{loc}}$.
\end{proposition}

\noindent{\bf Proof}\quad
Let $h\in C_c$ and let
$$(T_{\mu} h)(t)=t\int_{\mathbb R^+}h(t+s)\,d\mu(s)\qquad(t\in\mathbb R^+).$$
It follows as in \cite[Theorem~3.3.15]{Da:Book} or the proof of Proposition~\ref{pr:mol} that $T_{\mu} h$ is continuous. 
Also, we have $\text{supp}\,T_{\mu} h\subseteq\text{supp}\,h$, so $T_{\mu}$ maps $C_c$ into $C_c$. Moreover, for 
$h\in C_c$ and $\nu\in M_{\text{loc}}$ a calculation similar to the one in the proof of Proposition~\ref{pr:mol} shows that
$\langle h,\overline{D}_{\mu}\nu\rangle=\langle T_{\mu} h,\nu\rangle$.
Hence $\overline{D}_{\mu}=T_{\mu}^*$ and in particular $\overline{D}_{\mu}$ is weak-star continuous.
\ {\nopagebreak\hfill\raggedleft$\Box$\bigskip}

It follows from \cite[Proposition~8.4.30]{Pe-Bo} that the dual space of $L^1_{\text{loc}}$ can be identified with the inductive limit of the dual spaces $L^1[0,n)^*=L^{\infty}[0,n)$, which again can be identified with the space $L^{\infty}_c$ of measurable functions on $\mathbb R^+$ with compact support (with the inductive limit topology).

\begin{lemma}
\label{le:loc-weak}
The weak topology on $L^1_{\text{loc}}$ coincides with the topology $\tau$ on $L^1_{\text{loc}}$ obtained as the projective limit of the weak topologies on $L^1[0,n)$.
\end{lemma}

\noindent{\bf Proof}\quad
The topology $\tau$ is the coarsest topology on $L^1_{\text{loc}}$ making all the restrictions $R_n:L^1_{\text{loc}}\to(L^1[0,n),\text{weak})$ continuous. For  
$n\in\mathbb N$ and $\varphi\in L^{\infty}[0,n)=L^1[0,n)^*$, we let 
$U_{n,\varphi}=\{g\in L^1[0,n):|\langle g,\varphi\rangle|<1\}$. Then 
$\{U_{n,\varphi}:\varphi\in L^{\infty}[0,n)\}$ is a base for the 
$0$-neighbourhoods in the weak topology on $L^1[0,n)$. Hence 
$\{R_n^{-1}(U_{n,\varphi}):n\in\mathbb N,\ \varphi\in L^{\infty}[0,n)\}$
is a base for the $0$-neighbourhoods in the $\tau$ topology. Moreover, 
$R_n^{-1}(U_{n,\varphi})=\{f\in L^1_{\text{loc}}:|\langle R_nf,\varphi\rangle|<1\}
=\{f\in L^1_{\text{loc}}:|\langle f,S_n\varphi\rangle|<1\}$.
Since $(L^1_{\text{loc}})^*=L^{\infty}_c$ as the inductive limit of $L^{\infty}[0,n)$, it follows that $\tau$ equals the weak topology on $L^1_{\text{loc}}$.
{\nopagebreak\hfill\raggedleft$\Box$\bigskip}

It follows from the definition of the projective limit topology on $L^1_{\text{loc}}$ that a linear map $T:L^1_{\text{loc}}\to L^1_{\text{loc}}$ is continuous if and only if for every $n\in\mathbb N$ there exist $m\in\mathbb N$ and a constant $K$ such that $\|R_nTf\|\le K\|R_mf\|$ for every $f\in L^1_{\text{loc}}$. For weakly compact operators on $L^1_{\text{loc}}$ we have the following general description.

\begin{proposition}
\label{pr:loc-wcpt}
For a continuous operator $T$ on $L^1_{\text{loc}}$ the following conditions are equivalent:
\begin{enumerate}[(a)]
\item 
$T$ is weakly compact.
\item
There exists $m\in\mathbb N$ such that for every $n\in\mathbb N$ the operator $T_{nm}=R_nTS_m:L^1[0,m)\to L^1[0,n)$ is  weakly compact and moreover $Tf=0$ for every $f\in L^1_{\text{loc}}$ with $f=0$ on $[0,m)$.
\item 
There exists $m\in\mathbb N$ such that for every $n\in\mathbb N$ the operator $T_{nm}=R_nTS_m:L^1[0,m)\to L^1[0,n)$ is weakly compact and satisfies $R_nT=T_{nm}R_m$.
\end{enumerate}
\end{proposition}

\noindent{\bf Proof}\quad

(a)$\Rightarrow$(b):\quad
Clearly $T_{nm}$ is weakly compact for all $m,n\in\mathbb N$. Moreover, there exists a neighbourhood $U$ of 0 in $L^1_{\text{loc}}$ for which $T(U)$ is weakly relatively compact in $L^1_{\text{loc}}$. For every $n\in\mathbb N$ it follows that  $R_nT(U)$ is weakly relatively compact, in particular weakly bounded, and thus bounded in $L^1[0,n)$ by the principle of uniform boundedness.
There exists $m\in\mathbb N$ and $\delta>0$ such that $V=\{f\in L^1_{\text{loc}}:\|R_mf\|<\delta\}\subseteq U$. 
It follows that for every $n\in\mathbb N$ there exists a constant $K_n$ such that $\|R_nTf\|\le K_n\|R_mf\|$ for $f\in L^1_{\text{loc}}$. In particular, if $f\in L^1_{\text{loc}}$ with $f=0$ on $[0,m)$, then $Tf=0$ on $[0,n)$ for every $n\in\mathbb N$ and hence $Tf=0$. 

(b)$\Rightarrow$(c):\quad
Since $1-S_mR_m$ is the projection onto 
$\{f\in L^1_{\text{loc}}:f=0\text{ on }[0,m)\}$, we have $T(1-S_mR_m)=0$. Hence
$R_nT=R_nTS_mR_m=T_{nm}R_m$ for every $n\in\mathbb N$.

(c)$\Rightarrow$(a):\quad
Let $V=\{f\in L^1_{\text{loc}}:\|R_mf\|<1\}$. Then $V$ is a neighbourhood of 0 in $L^1_{\text{loc}}$ and $R_m(V)$ is the unit ball in $L^1[0,m)$, so $R_n(T(V))=T_{nm}(R_m(V))$ is weakly relatively compact in $L^1[0,n)$ for every $n\in\mathbb N$. By Lemma~\ref{le:loc-weak} and \cite[p.\,85]{Rob} it follows that $T(V)$ is weakly relatively compact in $L^1_{\text{loc}}$, so $T$ is weakly compact.
{\nopagebreak\hfill\raggedleft$\Box$\bigskip}

A slightly simpler version of the proof above shows the following result.

\begin{proposition}
\label{pr:loc-cpt}
For a continuous operator $T$ on $L^1_{\text{loc}}$ the following conditions are equivalent:
\begin{enumerate}[(a)]
\item 
$T$ is compact.
\item
There exists $m\in\mathbb N$ such that for every $n\in\mathbb N$ the operator $T_{nm}=R_nTS_m:L^1[0,m)\to L^1[0,n)$ is compact and moreover $Tf=0$ for every $f\in L^1_{\text{loc}}$ with $f=0$ on $[0,m)$.
\item 
There exists $m\in\mathbb N$ such that for every $n\in\mathbb N$ the operator $T_{nm}=R_nTS_m:L^1[0,m)\to L^1[0,n)$ is compact and satisfies $R_nT=T_{nm}R_m$.
\end{enumerate}
\end{proposition}

Since the derivations $D_{\mu}$ and the multipliers $T_{\mu}$ described in the beginning of the section are all 1-1 we obtain the next two results as consequences of Proposition~\ref{pr:loc-wcpt}.

\begin{corollary}
\label{co:loc-wcpt1}
There are no non-zero weakly compact derivations on $L^1_{\text{loc}}$.
\end{corollary}

\begin{corollary}
\label{co:loc-wcpt2}
There are no non-zero weakly compact multipliers on $L^1_{\text{loc}}$.
\end{corollary}

Motivated by Corollary~\ref{co:loc-wcpt1} we will now consider the weaker notions of (weakly) Montel derivations on $L^1_{\text{loc}}$ for which we have the following result similar to Theorem~\ref{th:cptvol}. 

\begin{theorem}
\label{th:monloc}
For $\mu\in M_{\text{loc}}$ the following conditions are equivalent:
\begin{enumerate}[(a)]
\item 
$D_{\mu}$ is a Montel derivation on $L^1_{\text{loc}}$.
\item 
$D_{\mu}$ is a weakly Montel derivation on $L^1_{\text{loc}}$.
\item
$\mu$ is absolutely continuous.
\item 
$\overline{D}_{\mu}$ is a Montel derivation on $M_{\text{loc}}$.
\item 
$\overline{D}_{\mu}$ is a weakly Montel derivation on $M_{\text{loc}}$.
\end{enumerate}
\end{theorem}

\noindent{\bf Proof}\quad
The implications (d)$\Rightarrow$(a)$\Rightarrow$(b) and (d)$\Rightarrow$(e)$\Rightarrow$(b) are obvious.

(b)$\Rightarrow$(c):\quad
Let $m\in\mathbb N$ and let $B_m$ be the closed unit ball in $L^1[0,m)$. Then 
$$S_m(B_m)=\{f\in L^1_{\text{loc}}:\text{supp}\,f\subseteq[0,m]\text{ and }\|R_mf\|\le1\},$$
so $R_n(S_m(B_m))\subseteq B_n$ for every $n\in\mathbb N$. Hence $S_m(B_m)$ is bounded in $L^1_{\text{loc}}$. It thus follows that $D_{\mu}(S_m(B_m))$ is weakly relatively compact in $L^1_{\text{loc}}$, so $R_m(D_{\mu}(S_m(B_m)))$ is weakly relatively compact in $L^1[0,m)$. Consequently $R_mD_{\mu} S_m$ is weakly compact. For $g\in L^1[0,m)$ and $t\in[0,m)$ we have
\begin{equation}
\label{eq:rds}
(R_mD_{\mu} S_m)(g)(t)=\int_{[0,t)}(t-s)(S_mg)(t-s)\,d\mu(s)
=\int_{[0,t)}(t-s)g(t-s)\,d(R_m\mu)(s).
\end{equation}
It is an easy corollary to the description of the continuous derivations on $L^1[0,1)$ mentioned in Section~\ref{sec:vol} and to Theorem~\ref{th:cptvol} that the continuous derivations on $L^1[0,m)$ are exactly the maps $\widetilde{D}_{\mu}g=Xg*\mu\ (g\in L^1[0,m))$ for some measure $\mu$ on $[0,m)$ with $t|\mu|([0,m-t))$ bounded as $t\to0_+$, and that $\widetilde{D}_{\mu}$ is (weakly) compact if and only if $\mu$ is absolutely continuous and $t|\mu|([0,m-t))\to0$ as $t\to0_+$. From \eqref{eq:rds} we see that
\begin{equation}
\label{eq:rds2}
R_mD_{\mu} S_m=\widetilde{D}_{R_m\mu},
\end{equation}
so we deduce that $R_m\mu$ is absolutely continuous. Since $m\in\mathbb N$ was arbitrary, this shows that $\mu$ is absolutely continuous on $\mathbb R^+$.

(c)$\Rightarrow$(d):\quad
Let $m\in\mathbb N$. It follows from the equivalent of \eqref{eq:rds2} for $\overline{D}_{\mu}$ and the comments in the proof of (b)$\Rightarrow$(c) that $R_m\overline{D}_{\mu} S_m$ is compact. For $\nu\in M_{\text{loc}}$ we observe that $R_m\overline{D}_{\mu}\nu$ only depends on $S_mR_m\nu$, that is $R_m\overline{D}_{\mu}\nu=R_m\overline{D}_{\mu} S_mR_m\nu$, so 
$R_m\overline{D}_{\mu}=R_m\overline{D}_{\mu} S_mR_m$.
Let $B$ be a bounded set in $M_{\text{loc}}$. Then $R_m(B)$ is bounded in $M[0,m)$, so
$$R_m(\overline{D}_{\mu}(B))=(R_m\overline{D}_{\mu} S_m)(R_m(B))$$
is relatively compact in $M[0,m)$. Since $m\in\mathbb N$ was arbitrary it follows from \cite[p.\,85]{Rob} that $\overline{D}_{\mu}(B)$ is relatively compact in $M_{\text{loc}}$. Hence $\overline{D}_{\mu}$ is Montel.
{\nopagebreak\hfill\raggedleft$\Box$\bigskip}

% {\bf [?? Note to self (others?): $R_mD_{\mu}$ compact for all $m\in\mathbb N$ does not imply that $D_{\mu}$ is compact, since the nbh depends on $m$.]}

\section{Derivations on $A(\omega)$}
\label{sec:aom}

% {\bf [?? Also properties of $(\omega_n)$; eg. $\omega_n(t)$ increasing or at least backwards shift continuous]}

In \cite{Pe:Fr} we studied the following class of weighted convolution Fr\'echet algebras (see \cite{Pe:Fr} for further details).
Let $\omega$ be an algebra weight on $\mathbb R^+$, that is, a positive Borel function satisfying: 
$\omega$ and $1/\omega$ are locally bounded on $\mathbb R^+$, 
$\omega$ is right continuous on $\mathbb R^+$,
$\omega$ is submultiplicative, that is $\omega(t+s)\le\omega(t)\omega(s)$ for $t,s\in\mathbb R^+$, and $\omega(0)=1$.
We then define $L^1({\omega})$ as the weighted space of functions $f$ on 
$\mathbb R^+$ for which $f\omega\in L^1(\mathbb R^+)$ with the norm 
$$\|f\|_{\omega}=\int_0^{\infty}|f(t)|\omega(t)\,dt.$$
It is well known that $L^1({\omega})$ with convolution product is a commutative Banach algebra. Similarly, we let $M(\omega)$ be the Banach algebra of locally finite complex Borel measures $\mu$ on $\mathbb R^+$ for which 
$$\|\mu\|_{\omega}=\int_{\mathbb R^+}\omega(t)\,d|\mu|t<\infty.$$
 
We consider an increasing sequence $\omega=(\omega_n)$ of algebra weights on $\mathbb R^+$ satisfying
\begin{enumerate}[(a)]
\item 
$\omega_n(t)\to\infty$ as $t\to\infty$ for every $n\in\mathbb N$,
\item
$\lim_{t\to\infty}\omega_n(t)^{1/t}=1$ for every $n\in\mathbb N$,
\item
$\sup_{t\in\mathbb R^+}\omega_{n+1}(t)/\omega_n(t)=\infty$ for every  
$n\in\mathbb N$.
\end{enumerate}
Let
$$A(\omega)=\bigcap_n L^1(\omega_n)\qquad\text{and}\qquad 
B(\omega)=\bigcap_nM(\omega_n)$$
and equip $A(\omega)$ and $B(\omega)$ with the increasing sequence of norms 
$\|\mu\|_n=\|\mu\|_{\omega_n}\ (\mu\in B(\omega))$. In this way $A(\omega)$ and $B(\omega)$ become Fr\'{e}chet algebras, which can be viewed as projective limits of $L^1(\omega_n)$ respectively $M(\omega_n)$.

In  \cite{Pe:Fr} we obtained the following characterisation of the derivations on $A(\omega)$.

\begin{theorem}[{\cite[Theorem~4.1]{Pe:Fr}}]
\label{th:der}
\ 
\begin{enumerate}[(a)]
\item 
Suppose that 
\begin{equation}
\label{eq:weco}
\text{for every $n\in\mathbb N$ there exists $m\in\mathbb N$ such that}\quad
\sup_{t\in\mathbb R^+}\frac{t\omega_n(t)}{\omega_m(t)}<\infty.
\end{equation}
Then 
$$D_{\mu}(f)=(Xf)\ast\mu\qquad(f\in A(\omega))$$
defines a continuous derivation on $A(\omega)$ for every $\mu\in B(\omega)$ and
conversely every derivation on $A(\omega)$ has this form. Also,
$\overline{D}_{\mu}(\nu)=(X\nu)\ast\mu$ for $\nu\in B(\omega)$ extends $D_{\mu}$ to a continuous derivation on $B(\omega)$.
\item
If condition \eqref{eq:weco} is not satisfied, then there are no
non-zero derivations on $A(\omega)$.
\end{enumerate}
\end{theorem}

(In particular, derivations on $A(\omega)$ are automatically continuous.) 

As described in \cite[Proposition~3.3]{Pe:Fr} the algebra $B(\omega)$ can be identified with the dual space of the space 
$D(1/\omega)=\cup_{n\in\mathbb N}C_0(1/\omega_n)$ with the inductive limit topology
(where $C_0(1/\omega_n)$ is space of continuous functions $h$ on
$\mathbb R^+$ for which $h/\omega_n$ is vanishes at infinity).
We will first show that the derivations $\overline{D}_{\mu}$ are weak-star continuous.

\begin{proposition}
\label{pr:aom}
Assume that condition \eqref{eq:weco} is satisfied and let $\mu\in B(\omega)$. Then the derivation $\overline{D}_{\mu}$ is weak-star continuous on $B(\omega)$.
\end{proposition}

\noindent{\bf Proof}\quad
Let $h\in D(1/\omega)$ and let
$$(T_{\mu} h)(t)=t\int_{\mathbb R^+}h(t+s)\,d\mu(s)\qquad(t\in\mathbb R^+).$$
As in the proofs of Propositions~\ref{pr:mol} and \ref{pr:locwks} it follows that
$T_{\mu} h$ is continuous on $\mathbb R^+$. Choose $n\in\mathbb N$ such that $h\in C_0(1/\omega_n)$, and then choose 
$m\in\mathbb N$ and $C>0$ such that $t\omega_n(t)\le C\omega_m(t)$ for $t\in\mathbb R^+$. Let $\varepsilon$ be a decreasing function with $\varepsilon(t)\to0$ as $t\to\infty$ such that $|h(t)|\le\varepsilon(t)\omega_n(t)$ for $t\in\mathbb R^+$. We then have
\begin{align*}
|(T_{\mu} h)(t)|
&\le t\int_{\mathbb R^+}\varepsilon(t+s)\omega_n(t+s)\,d|\mu|(s)\\
&\le \varepsilon(t)t\omega_n(t)\int_{\mathbb R^+}\omega_n(s)\,d|\mu|(s)\\
&\le C\varepsilon(t)\omega_m(t)\|\mu\|_n
\end{align*}
for $t\in\mathbb R^+$. 
Hence $T_{\mu} h\in C_0(1/\omega_m)$, and we deduce that $T_{\mu}$ is a continuous linear operator on $D(1/\omega)$. Moreover, for $\nu\in B(\omega)$ and 
$h\in D(1/\omega)$ a calculation similar to the one in the proof of Proposition~\ref{pr:mol} shows that
$\langle h,\overline{D}_{\mu}\nu\rangle=\langle T_{\mu} h,\nu\rangle$.
Hence $\overline{D}_{\mu}=T_{\mu}^*$ and in particular $\overline{D}_{\mu}$ is weak-star continuous.
{\nopagebreak\hfill\raggedleft$\Box$\bigskip}

An argument similar to the one above shows that for a derivation $D_{\mu}$ on $L^1(\omega)$, the extension $\overline{D}_{\mu}$ is weak-star continuous on $M(\omega)$.

As for $L^1_{\text{loc}}$ in the previous section we will show that the zero operator is the only weakly compact derivation on $A(\omega)$. 

\begin{theorem}
\label{th:aom}
There are no non-zero weakly compact derivations on $A(\omega)$.
\end{theorem}

\noindent{\bf Proof}\quad
We have $A(\omega)^*=\bigcup_n L^{\infty}(1/\omega_n)$ by \cite[Corollary~2.2]{Pe:Fr}, and it follows as in the proof of Lemma~\ref{le:loc-weak} that the weak topology on $A(\omega)$ coincides with the topology $\tau$ obtained as the projective limit of the weak topologies on $L^1(\omega_n)$. 
% The topology $\tau$ is the coarsest topology on $A(\omega)$ making all the inclusions 
% $\iota_n:A(\omega)\to(L^1(\omega_n),\text{weak})$ continuous. For  
% $n\in\mathbb N$ and $\varphi\in L^{\infty}(1/\omega_n)=L^1(\omega_n)^*$, we let 
% $U_{n,\varphi}=\{f\in L^1(\omega_n):|\langle f,\varphi\rangle|<1\}$. Then 
% $\{U_{n,\varphi}:\varphi\in L^{\infty}(1/\omega_n)\}$ is a base for the 
% $0$-neighbourhoods in the weak topology on $L^1(\omega_n)$. Hence 
% $\{\iota_n^{-1}(U_{n,\varphi}):n\in\mathbb N,\ \varphi\in L^{\infty}(1/\omega_n)\}$
% is a base for the $0$-neighbourhoods in the $\tau$ topology. Moreover, 
% $\iota_n^{-1}(U_{n,\varphi})=\{f\in A(\omega):|\langle f,\varphi\rangle|<1\}$.
%  that $\tau$ equals the weak topology on $A(\omega)$.
%\medskip

Let $D$ be a weakly compact derivation on $A(\omega)$. There exists a neighbourhood $U$ in $A(\omega)$ for which $D(U)$ is weakly relatively compact in $A(\omega)$. By the above it follows that $D(U)$ is weakly relatively compact in $L^1(\omega_n)$ for every $n\in\mathbb N$. In particular $D(U)$ is weakly bounded and by the principle of uniform boundedness thus bounded in $L^1(\omega_n)$ for every $n\in\mathbb N$. There exists $m\in\mathbb N$ and $\delta>0$ such that 
$V=\{f\in A(\omega):\|f\|_{L^1(\omega_m)}<\delta\}\subseteq U$. Let 
$n\in\mathbb N$. It follows that there exists a constant $K_n$ such that 
$\|Df\|_{L^1(\omega_n)}\le K_n\|f\|_{L^1(\omega_m)}$ for $f\in A(\omega)$. Since 
$A(\omega)$ is dense in $L^1(\omega_n)$ we deduce that $D$ extends to a continuous linear operator $D_n:L^1(\omega_m)\to L^1(\omega_n)$. In particular $D_m$ is a derivation on $L^1(\omega_m)$, so by Johnson's result
(\cite{Jo:Cont} or \cite[Theorem~5.2.32]{Da:Book}) we have $D_m=0$ and thus $D=0$.
{\nopagebreak\hfill\raggedleft$\Box$\bigskip}

We finish the paper by showing that under a slightly stronger assumption than \eqref{eq:weco}, the Montel derivations $D_{\mu}$ on $A(\omega)$ correspond to absolutely continuous measures $\mu$ 
(as for $L^1_{\text{loc}}$).

\begin{theorem}
\label{th:monaom}
Suppose that
$$\text{for every $n\in\mathbb N$ there exists $m\in\mathbb N$ such that}\quad
\frac{t\omega_n(t)}{\omega_m(t)}\to0\text{ as }t\to\infty$$
and let $\mu\in B(\omega)$. Then the following conditions are equivalent:
\begin{enumerate}[(a)]
\item 
$D_{\mu}$ is a Montel derivation on $A(\omega)$.
\item 
$D_{\mu}$ is a weakly Montel derivation on $A(\omega)$.
\item
$\mu$ is absolutely continuous.
\item 
$\overline{D}_{\mu}$ is a Montel derivation on $B(\omega)$.
\item 
$\overline{D}_{\mu}$ is a weakly Montel derivation on $B(\omega)$.
\end{enumerate}
\end{theorem}

\noindent{\bf Proof}\quad
The implications (d)$\Rightarrow$(a)$\Rightarrow$(b) and (d)$\Rightarrow$(e)$\Rightarrow$(b) are obvious.

(b)$\Rightarrow$(c):\quad
Let $(e_k)$ be a bounded approximate identity for $A(\omega)$ and let $t>0$. Then
$\{\delta_t*e_k:k\in\mathbb N\}$ is bounded in $A(\omega)$, so there exist a subsequence $(e_{k_j})$ and $f\in A(\omega)$ such that
$$D_{\mu}(\delta_t*e_{k_j})\to f\qquad\text{weakly in $A(\omega)$ as $j\to\infty$}.$$
Now, proceed as in the proof of Theorem~\ref{th:cptvol}\,(b)$\Rightarrow$(c).

(c)$\Rightarrow$(d):\quad
Let $E$ be a bounded set in $B(\omega)$. 
Let $n\in\mathbb N$ and choose $m\in\mathbb N$ such that $\overline{D}_{\mu}$ extends to a continuous linear map 
$\overline{D}_{mn}:M(\omega_m)\to M(\omega_n)$. We may assume that $m\in\mathbb N$ is chosen so that $t\omega_n(t)/\omega_m(t)\to0$ as $t\to\infty$. 
It follows from the proof of \cite[Theorem~4.1\,(c)$\Rightarrow$(d)]{De-Gh-Gr} 
(see also the proof of Theorem~\ref{th:cptvol}) combined with the estimate
$$\left\|\frac{\overline{D}_{mn}\delta_t}{\omega_m(t)}\right\|_n
=\left\|\frac{t}{\omega_m(t)}\,\delta_t*\mu\right\|_n
\le\left\|\frac{t}{\omega_m(t)}\,\delta_t\right\|_n\cdot\left\|\mu\right\|_n
=\frac{t\omega_n(t)}{\omega_m(t)}\cdot\left\|\mu\right\|_n\to0$$
as $t\to\infty$ that $\overline{D}_{mn}$ maps the unit ball in $M(\omega_m)$ to a relatively compact set in $M(\omega_n)$. 
Since $E$ is bounded in $M(\omega_m)$, we deduce that $\overline{D}_{\mu}(E)=\overline{D}_{mn}(E)$ is  relatively compact in $M(\omega_n)$. Finally, by \cite[p.\,85]{Rob} this proves that $\overline{D}_{\mu}(E)$ is relatively compact in $B(\omega)$, so $\overline{D}_{\mu}$ is Montel.
{\nopagebreak\hfill\raggedleft$\Box$\bigskip}

\bigskip
\bigskip

% \bibliography{myref}

% \bibliographystyle{plain}

\bigskip

\noindent
Thomas Vils Pedersen\\
Department of Mathematical Sciences\\
University of Copenhagen\\
Universitetsparken 5\\
DK-2100 Copenhagen \O\\
Denmark\\
vils@math.ku.dk

\end{document}